\documentclass[reqno]{amsart}
\usepackage{amssymb}
\usepackage{epsf}

\newtheorem{theorem}{Theorem}

\theoremstyle{definition}

\theoremstyle{remark}
\newtheorem{remark}{Remark}

\newcommand{\N}{\mathbb N}

\newcommand{\C}{\mathbb C}
\newcommand{\R}{\mathbb R}
\newcommand{\Z}{\mathbb Z}

\newcommand{\T}{\mathbb T}
\renewcommand{\L}{\mathbb L}
\renewcommand{\P}{\mathcal P}

\begin{document}

\title[Elliptic Beta Integrals]
{Short proofs of the elliptic beta integrals}

\author{V.P. Spiridonov}

\address{Bogoliubov Laboratory of Theoretical Physics, JINR,
Dubna, Moscow reg. 141980, Russia}

\thanks{{\em Date:} March 2004; to be published in the {\em Ramanujan J.} \\
 This work is supported in part by the Russian Foundation for Basic
Research (RFBR) grant no. 03-01-00781}

\begin{abstract}
We give elementary proofs of the univariate elliptic beta
integral with bases $|q|, |p|<1$ and its multiparameter generalizations
to integrals on the $A_n$ and $C_n$ root systems.
We prove also some new unit circle multiple elliptic beta integrals,
which are well defined for $|q|=1$, and their $p\to 0$ degenerations.
\end{abstract}

\dedicatory{Dedicated to Richard Askey on the occasion of his seventieth
birthday}

\maketitle

\tableofcontents

\section{Introduction}

Hypergeometric type functions play a major role in the theory
of special functions, see, e.g., \cite{aar:special}. The concept
of elliptic hypergeometric integrals was introduced by the author
in \cite{spi:elliptic}, where a univariate elliptic
analog of the beta type integrals has been discovered. First multidimensional
extensions of this integral were proposed by van Diejen and the author
in \cite{die-spi:elliptic} (two types of integrals for the $C_n$ root system).
Other generalizations were introduced in \cite{spi:elliptic} (another $C_n$ integral
and three different integrals for the $A_n$ root system) and \cite{spi-war:inversions}
(a new $A_n$ integral which, being related to the series prescribed earlier to the
$D_n$ root system, can also be considered formally as a $D_n$ integral).
These are all elliptic beta integrals found so far.
The first complete proofs of the multiparameter $C_n$
integral of \cite{die-spi:elliptic} and $A_n$ integral of \cite{spi:elliptic}
were obtained recently by Rains \cite{rai:trans}. The present paper describes
elementary proofs of these exact integration formulas modelled along the
proof of the Askey-Wilson integral \cite{ask-wil:some}
given by Wilf and Zeilberger in \cite{wz:invent}.
In addition, we construct some new multiparameter elliptic beta
integrals well defined for $|q|=1$ and prove them together with
their $q$-hypergeometric degenerations. We would like to mention also that
the $A_n$ (or $``D_n"$) integral introduced by Warnaar and the author in
\cite{spi-war:inversions} has been proved there using similar
elementary means.

We take two base variables $q, p\in\mathbb{C}$ satisfying
constraints $|q|,|p|<1$. The key Jacobi type theta function has
the form
$$
\theta(z;p)=(z;p)_\infty(pz^{-1};p)_\infty,
$$
where $(a;p)_\infty=\prod_{j=0}^\infty(1-ap^j)$.
Its main transformation properties are
$$
\theta(z^{-1};p)=\theta(pz;p)=-z^{-1}\theta(z;p).
$$
The standard elliptic gamma function \cite{rui:first}
\begin{equation}
\Gamma(z;q,p) = \prod_{j,k=0}^\infty\frac{1-z^{-1}q^{j+1}p^{k+1}}
{1-zq^jp^k}
\label{ell-gamma}\end{equation}
is symmetric in bases $q,p$ and satisfies equations
$$
\Gamma(qz;q,p)=\theta(z;p)\Gamma(z;q,p),\qquad
\Gamma(pz;q,p)=\theta(z;q)\Gamma(z;q,p)
$$
together with the reflection equation $\Gamma(z;q,p)\Gamma(pq/z;q,p)=1$.
We shall drop bases $q$ and $p$ from the notation $\Gamma(z;q,p)$ and use conventions
\begin{eqnarray*}
&& \Gamma(tz^\pm)=\Gamma(tz,tz^{-1}),\quad
\Gamma(z^{\pm 2})=\Gamma(z^2,z^{-2}),
\\
&& \Gamma(t_1,\ldots,t_m)\equiv\prod_{r=1}^m \Gamma(t_r;q,p).
\end{eqnarray*}
Similarly, for theta functions we assume that
$$
\theta(t_1,\ldots,t_m;p)=\theta(t_1;p)\cdots\theta(t_m;p),
\quad \theta(tz^\pm;p)=\theta(tz,tz^{-1};p).
$$

It is convenient to introduce the exponential parametrization:
\begin{eqnarray}\nonumber
&& q=e^{2\pi i\frac{\omega_1}{\omega_2}}, \qquad p=e^{2\pi
i\frac{\omega_3}{\omega_2}},\qquad r=e^{2\pi
i\frac{\omega_3}{\omega_1}},
\\
&& \tilde q =e^{-2\pi i\frac{\omega_2}{\omega_1}}, \qquad \tilde
p=e^{-2\pi i\frac{\omega_2}{\omega_3}},\qquad \tilde r=e^{-2\pi
i\frac{\omega_1}{\omega_3}},
\label{ell-bases}\end{eqnarray}
where $\omega_{1,2,3}$ are some complex numbers.
For $|q|, |p|, |r|<1$, the modified elliptic gamma function has the form
\cite{spi:elliptic}
\begin{equation}
G(u;\boldsymbol{\omega})= \prod_{j,k=0}^\infty \frac{(1-e^{-2\pi
i\frac{u}{\omega_2}}q^{j+1}p^{k+1}) (1-e^{2\pi i
\frac{u}{\omega_1}}{\tilde q}^{j+1}{r }^k)} {(1-e^{2\pi i
\frac{u}{\omega_2}}q^jp^k) (1-e^{-2\pi i
\frac{u}{\omega_1}}{\tilde q}^j{r}^{k+1})}
\label{ell-d}\end{equation}
and satisfies difference equations
\begin{eqnarray}
&& G(u+\omega_1;\boldsymbol{\omega})=\theta(e^{2\pi
i\frac{u}{\omega_2}};p) G(u;\boldsymbol{\omega}),
\label{ell-1eq} \\
&& G(u+\omega_2;\boldsymbol{\omega})=\theta(e^{2\pi
i\frac{u}{\omega_1}};r) G(u;\boldsymbol{\omega}),
\label{ell-2eq} \\
&& G(u+\omega_3;\boldsymbol{\omega})=e^{-\pi iB_{2,2}(u;\mathbf{\omega})}
G(u;\boldsymbol{\omega}),
\label{ell-3eq-r}\end{eqnarray}
where
$$
B_{2,2}(u;\boldsymbol{\omega})=\frac{u^2}{\omega_1\omega_2}
-\frac{u}{\omega_1}-\frac{u}{\omega_2}+
\frac{\omega_1}{6\omega_2}+\frac{\omega_2}{6\omega_1}+\frac{1}{2}.
$$

The representation \cite{die-spi:elliptic}
\begin{equation}
G(u;\boldsymbol{\omega})=e^{-\pi iP(u)}\Gamma(e^{-2\pi
i\frac{u}{\omega_3}}; \tilde r, \tilde p),
\label{gamma-tr}\end{equation} where $P(u)$ is the following
polynomial of the third degree
\begin{equation}
P\Big(u+\sum_{k=1}^3\frac{\omega_k}{2}\Big)=\frac{u(u^2-\frac{1}{4}\sum_{k=1}^3\omega_k^2)}
{3\omega_1\omega_2\omega_3},
\label{p(u)}\end{equation}
shows that the modified elliptic gamma function is well defined
for $|q|=1$ with $\omega_1/\omega_2>0$
in distinction from $\Gamma(z;q,p)$. It is related to modular transformations
for the standard elliptic gamma function \cite{fel-var:elliptic}.

We have the symmetry $G(u;\omega_1,\omega_2,\omega_3)=G(u;\omega_2,\omega_1,\omega_3)$
and the reflection equation
$$
G(a;\boldsymbol{\omega})G(b;\boldsymbol{\omega})=1,\quad a+b=\sum_{k=1}^3\omega_k,
$$
following from the property $P(\sum_{k=1}^3\omega_k-u)=-P(u)$.
In the limit $p,r\to 0$, the function $G(u;\boldsymbol{\omega})$
becomes reciprocal to the double sine function
\begin{equation}
\lim_{p,r\to 0} \frac{1}{G(u;\boldsymbol{\omega})}
=S(u;\omega_1,\omega_2) = \frac{(e^{2\pi i u/\omega_2}; q)_\infty}
{(e^{2\pi iu/\omega_1}\tilde q; \tilde q)_\infty}.
\label{2d-sin}\end{equation}

We use also the notation
\begin{eqnarray*}
&&
G(u_1,\ldots,u_m)\equiv G(u_1;\boldsymbol{\omega})\cdots
G(u_1;\boldsymbol{\omega}),
\\ &&
S(u_1,\ldots,u_m)\equiv S(u_1;\boldsymbol{\omega})\cdots
S(u_1;\boldsymbol{\omega}).
\end{eqnarray*}

For $\text{Im}(\omega_{1}/\omega_2)>0$ (or $|q|<1$) the
double sine function  can be defined
as the infinite product (\ref{2d-sin}). However, it remains a well
defined meromorphic function of $u$ in the domain $\omega_1/\omega_2>0$
(so that $|q|=1$), see, e.g., \cite{kls:unitary}. Its zeros are located at
the points $u=-\omega_1\mathbb{N}-\omega_2\mathbb{N}$ and poles occupy the
lattice $u=\omega_1(1+\mathbb{N})+\omega_2 (1+\mathbb{N})$.
Asymptotically,
\begin{eqnarray}\label{asymp1}
&& \lim_{\text{Im}(\frac{u}{\omega_1}),\text{Im}(\frac{u}{\omega_2})\to
+\infty}S(u;\boldsymbol{\omega})=1,
\\ &&
\lim_{\text{Im}(\frac{u}{\omega_1}),\text{Im}(\frac{u}{\omega_2})\to
-\infty} e^{\pi iB_{2,2}(u;\boldsymbol{\omega})}
S(u;\boldsymbol{\omega}) =1 .
\label{asymp2}\end{eqnarray}

\section{The univariate elliptic beta integral}

We take $z\in\C$ and five complex parameters $t_m, m=1, \dots,5,$
assume that $|q|,|p|<1$ and compose the kernel function
\begin{equation}
\rho(z,t_1,\ldots,t_5)=\frac{\prod_{m=1}^5\Gamma(t_mz^\pm,At_m^{-1})}
{\Gamma(z^{\pm 2},Az^\pm)\prod_{1\leq m<s\leq 5}\Gamma(t_mt_s)},
\label{kernel}\end{equation}
where $A=\prod_{m=1}^5t_m$. This function has sequences of poles
which converge to zero along the points
$$
\P=\{t_mq^jp^k,\, A^{-1}q^{j+1}p^{k+1}
\}_{m=1,\ldots, 5,\, j,k \in\N}
$$
and diverge to infinity along their reciprocals $\P^{-1}$.
Denote as $C$ a contour on the complex plane with positive orientation
which separates $\P$ and $\P^{-1}$. For instance, for
$|t_m|<1,\; |pq|<|A|,$
the contour $C$ may coincide with the unit circle $\mathbb{T}$.

\begin{theorem}
\begin{equation}\label{ell-int}
\int_C\rho(z,t_1,\ldots,t_5)\frac{dz}{z}=\frac{4\pi i}
{(q;q)_\infty(p;p)_\infty}.
\end{equation}
\end{theorem}
\begin{proof}
The first step consists in establishing the following
$q$-difference equation for the kernel function:
\begin{equation}
\rho(z,qt_1,t_2,\ldots,t_5)-\rho(z,t_1,\ldots,t_5)
=g(q^{-1}z,t_1,\ldots,t_5)-g(z,t_1,\ldots,t_5),
\label{eqn}\end{equation}
where
\begin{equation}
g(z,t_1,\ldots,t_5)=\rho(z,t_1,\ldots,t_5)
\frac{\prod_{j=1}^5\theta(t_jz;p)}{\prod_{j=2}^5\theta(t_1t_j;p)}
\frac{\theta(t_1A;p)}{\theta(z^2,Az;p)}\frac{t_1}{z}.
\label{g}\end{equation}
In the limits $p\to 0$ and subsequent $t_5\to 0$, this equation is reduced
to that used in \cite{wz:invent} for proving the Askey-Wilson integral
\cite{ask-wil:some}. The simple $p\to 0$ limit reduces equality (\ref{ell-int})
to the $q$-beta integral introduced by Rahman in \cite{rah:integral}.
After division of equation (\ref{eqn})
by $\rho(z,t_1,\ldots,t_5)$ it takes the form
\begin{eqnarray}\nonumber
\lefteqn{\frac{\theta(t_1z,t_1z^{-1};p)}{\theta(Az,Az^{-1};p)}
\prod_{j=2}^5\frac{\theta(At_j^{-1};p)}{\theta(t_1t_j;p)}-1 } &&
\\ &&
=\frac{t_1\theta(t_1A;p)}{z\theta(z^2;p)\prod_{j=2}^5\theta(t_1t_j;p)}
\left(\frac{z^4\prod_{j=1}^5\theta(t_jz^{-1};p)}{\theta(Az^{-1};p)} -
\frac{\prod_{j=1}^5\theta(t_jz;p)}{\theta(Az;p)}\right).
\label{eqn-exp}\end{eqnarray}

Both sides of this equality represent elliptic functions of $\log z$
(i.e., they are invariant under the transformation $z\to pz$) with
equal sets of poles and their residues, e.g.,
$$
\lim_{z\to A}\theta(Az^{-1};p)\, (\text{l.h.s.}) =
\frac{\theta(t_1A,t_1A^{-1};p)}{\theta(A^2;p)}
\prod_{j=2}^5\frac{\theta(At_j^{-1};p)}{\theta(t_1t_j;p)}
$$
with the same result for the right-hand side.
Therefore, the difference of the expressions on two sides of (\ref{eqn-exp})
is an elliptic function without poles, that is a constant. The latter
constant is equal to zero since for $z=t_1$ validity of (\ref{eqn-exp})
is evident.

Now we integrate equation (\ref{eqn}) over the variable $z\in C$
and obtain
\begin{equation}
I(qt_1,t_2,\ldots,t_5)-I(t_1,\ldots,t_5)=
\left(\int_{q^{-1}C}-\int_C\right)g(z,t_1,\ldots,t_5)
\frac{dz}{z},
\label{int-eqn}\end{equation}
where $I(t_1,\ldots,t_5)=\int_C\rho(z,t_1,\ldots,t_5)dz/z$ and
$q^{-1}C$ denotes the contour $C$ scaled by $q^{-1}$ with respect
to the $z=0$ point. Function (\ref{g}) has converging sequences of
poles at $z=\{t_mq^jp^k, A^{-1}q^jp^{k+1}\}$
and diverging ones at $z=\{ t_m^{-1}q^{-j-1}p^{-k},$
$Aq^{-j-1}p^{-k-1}\}$ for $m=1,\ldots,5$ and $j,k\in\mathbb{N}$.
Taking $C=\mathbb{T}$, we see that for $|t_m|<1$ and $|p|<|A|$ there
are no poles in the annulus $1\leq |z|\leq |q|^{-1}$. Therefore,
we can deform $q^{-1}\mathbb{T}$ to $\mathbb{T}$
in (\ref{int-eqn}) and get zero on the right-hand side
yielding $I(qt_1,t_2,\ldots,t_5)=I(t_1,\ldots,t_5)$.
Assuming that $|p|,|q|<|A|$, we  obtain by symmetry that
$I(pt_1,t_2,\ldots,t_5)=$ $I(t_1,\ldots,t_5)$. Further transformations
$t_1\to q^{\pm1}t_1$ and $t_1\to p^{\pm1}t_1$ can be
performed only if they keep parameters inside the annulus
of analyticity of the function $I(t_1,\ldots,t_5)$.

Temporarily, we take real $p$ and $q$, $p<q$, such that $p^n\neq q^k$ for any
$n,k\in\mathbb{N}$ and assume  that the arguments of $t_m^{\pm1},\, m=1,\ldots,5,$
and $A^{\pm 1}$ differ pairwise. Now we fix $C$ to be a contour that
encircles $\P$ and two cuts $c_1=[t_1, t_1p^2]$, $c_2=[(pq/A)p^{-2},pq/A]$
and excludes their reciprocals. Now we scale $t_1\to t_1q^k$, $k=1,2,\ldots,$
and, as soon as $t_1q^k$ enters the interval $[t_1p,t_1p^2]$, we perform the
scaling $t_1\to t_1p^{-1}$ which does not take parameters outside $c_1$ or $c_2$.
In this way we obtain $I(q^jp^{-k}t_1,t_2,\ldots,t_5)=I(t_1,\ldots,t_5)$
for all $j,k\in\N$ such that $q^jp^{-k}\in [1,p]$. Since such a set of points
is dense, we conclude that $I$ does not depend on $t_1$ and, by symmetry, on
all $t_m$.

Alternatively, we can expand $I(t_1,\ldots,t_5)=$
$\sum_{n=0}^\infty I_n(t_1,\ldots,t_5)p^n$
and find that $I_n(qt_1,\ldots,t_5)=I_n(t_1,\ldots,t_5)$ termwise, similar to
the situations described earlier in \cite{die-spi:elliptic,spi:elliptic}.
The coefficients $I_n$ are analytic in the parameters near the $t_m=0$ points
(it is the convergence of the $p$-expansion that imposes constraints on the
absolute values of parameters from below). Therefore we can iterate $t_1\to qt_1$
scalings until reaching the limiting point. As a result, the coefficients
$I_n$ and the integral $I$ itself do not depend on $t_1$ and, consequently,
on all $t_m$.

Thus $I$ is a constant depending only on $p$
and $q$. Its value, given by the right-hand side of (\ref{ell-int}),
is found after forcing $C$ to cross over the poles at $z=t_1^{\pm 1}$, picking
up the corresponding residues, and taking the limit $t_2\to 1/t_1$,
like it was done in the residue calculus of \cite{die-spi:elliptic}.
After proving the integration formula for a restricted region of parameters, we
can continue it analytically to the parameters domain allowed by the contour $C$.
\end{proof}

\begin{remark}
Integral (\ref{ell-int}) was proven in \cite{spi:elliptic} with the help of the
Bailey's $_2\psi_2$ sum and an analytical continuation in a discrete
set of parameters. Here we gave an elementary proof which did
not use any $q$-series or integral identities and used
analytical continuations in a minimal way.
\end{remark}

\section{A multiparameter $C_n$ integral}

In order to treat the multiparameter elliptic beta integral for the
$C_n$ root system, which was introduced by van Diejen and the author in
\cite{die-spi:elliptic} and tagged there as a type I integral, we take
$z=(z_1,\ldots,z_n)\in\C^n$ and $2n+3$ complex parameters
$t=(t_1,\ldots,t_{2n+3})$ and compose the kernel function
\begin{equation}
\rho(z,t;C_n)=\prod_{1\leq i<j\leq n}\frac{1}{\Gamma(z_i^\pm z_j^\pm)}
\prod_{i=1}^n\frac{\prod_{m=1}^{2n+3}\Gamma(t_mz_i^\pm)}
{\Gamma(z_i^{\pm 2},Az_i^\pm)}\frac{\prod_{m=1}^{2n+3}
\Gamma(At_m^{-1})}{\prod_{1\leq m<s\leq 2n+3}\Gamma(t_mt_s)},
\label{kernel-C}\end{equation}
where $A=\prod_{m=1}^{2n+3}t_m$. We denote
$$
\P=\{t_mq^jp^k,\, A^{-1}q^{j+1}p^{k+1}
\}_{m=1,\ldots, 2n+3,\, j,k \in\N}
$$
the set of points on the complex plane along which the
poles of (\ref{kernel-C}) in $z_i$ converge to zero; $C$ is a positively oriented
contour separating $\P$ and $\P^{-1}$; $dz/z=\prod_{i=1}^ndz_i/z_i$.

\begin{theorem}
\begin{equation}\label{ell-int-C}
\int_{C^n}\rho(z,t;C_n)\frac{dz}{z}=\frac{2^nn!(2\pi i)^n}
{(q;q)_\infty^n(p;p)_\infty^n}.
\end{equation}
\end{theorem}
\begin{proof}
The kernel function satisfies a $q$-difference equation, analogous to (\ref{eqn}),
\begin{eqnarray} \nonumber
&& \rho(z,qt_1,t_2,\ldots,t_{2n+3};C_n)-\rho(z,t;C_n)
\\ && \makebox[4em]{}
=\sum_{i=1}^n\left(g_i(z_1,...,q^{-1}z_i,\ldots,z_n,t)-g_i(z,t)\right),
\label{eqn-C}\end{eqnarray}
where
\begin{equation}
g_i(z,t)=\rho(z,t;C_n)\prod_{j=1,\neq i}^n\frac{\theta(t_1z_j^\pm;p)}
{\theta(z_iz_j^\pm;p)}\frac{\prod_{j=1}^{2n+3}\theta(t_jz_i;p)}
{\prod_{j=2}^{2n+3}\theta(t_1t_j;p)}
\frac{\theta(t_1A;p)}{\theta(z_i^2,Az_i;p)}\frac{t_1}{z_i}.
\label{g-C}\end{equation}
Dividing (\ref{eqn-C}) by $\rho(z,t;C_n)$, we obtain
\begin{eqnarray}\nonumber
\lefteqn{\prod_{i=1}^n\frac{\theta(t_1z_i^\pm;p)}{\theta(Az_i^\pm;p)}
\prod_{j=2}^{2n+3}\frac{\theta(At_j^{-1};p)}{\theta(t_1t_j;p)}-1
=\frac{t_1\theta(t_1A;p)}{\prod_{j=2}^{2n+3}\theta(t_1t_j;p)}
\sum_{i=1}^n\frac{1}{z_i\theta(z_i^2;p)}
}&&
\\ &&
\times
\prod_{j=1,\neq i}^n\frac{\theta(t_1z_j^\pm;p)}{\theta(z_iz_j^\pm;p)}
\left(z_i^{2n+2}\frac{\prod_{j=1}^{2n+3}\theta(t_jz_i^{-1};p)}
{\theta(Az_i^{-1};p)} -
\frac{\prod_{j=1}^{2n+3}\theta(t_jz_i;p)}
{\theta(Az_i;p)}\right).
\label{eqn-exp-C}\end{eqnarray}
Both sides of this equality are invariant under the transformation $z_1\to pz_1$
and have equal sets of poles
(the poles at $z_1=z_j,z_j^{-1}, j=2,\ldots, n,$ and $z_1=\pm p^{k/2},\ k\in\N$
on the right-hand side cancel each other) and their residues, e.g.,
$$
\lim_{z_1\to A}\theta(Az_1^{-1};p)\, (\text{l. or r. h.s.}) =
\frac{\theta(t_1A^\pm;p)}{\theta(A^2;p)}\prod_{j=1,\neq i}^n\frac{
\theta(t_1z_i^\pm;p)}{\theta(Az_i^\pm;p)}
\prod_{j=2}^{2n+3}\frac{\theta(At_j^{-1};p)}{\theta(t_1t_j;p)}.
$$
Therefore, the functions on two sides of (\ref{eqn-exp-C})
differ only by a constant independent on $z_1$ which is equal to zero
since for $z_1=t_1$ equality (\ref{eqn-exp-C}) is evident.

Integrating equality (\ref{eqn-C}) over the variables $z\in C^n$, we obtain
\begin{equation}
I(qt_1,t_2,\ldots,t_{2n+3})-I(t)=\sum_{i=1}^n
\left(\int_{C^{i-1}\times(q^{-1}C)\times C^{n-i}}-\int_{C^n}\right)g_i(z,t)\frac{dz}{z},
\label{int-eqn-C}\end{equation}
where $I(t)=\int_{C^n}\rho(z,t;C_n)dz/z$ and
$q^{-1}C$ is the scaled partner of $C$.

Poles of functions (\ref{g-C}) in $z_i$ converge to zero along the set $z_i=\{t_mq^jp^k,$
$A^{-1}q^jp^{k+1}\}$ and diverge
for $z_i=\{t_m^{-1}q^{-j-1}p^{-k},Aq^{-j-1}p^{-k-1}\}$, where $m=1,\ldots,2n+3$,
$j,k\in\mathbb{N}$. For $|t_m|<1$ and $|p|<|A|$ the region
$1\leq |z_i|\leq |q|^{-1}$ does not contain poles, so that we can
take $C=\T$, deform $q^{-1}\mathbb{T}$ to $\mathbb{T}$
in (\ref{int-eqn-C}) and obtain $I(qt_1,t_2,\ldots,t_{2n+3})=I(t)$.
Repeating almost literally the analytical continuation procedure
used in the $n=1$ case, we see that $I$ is a function of $p$ and $q$ only.
The latter is found from the residue calculus of \cite{die-spi:elliptic}
to be equal to the right-hand side of (\ref{ell-int-C}).
\end{proof}

\begin{remark}
Integration formula (\ref{ell-int-C}) was conjectured in \cite{die-spi:elliptic}
with many justifying arguments (e.g., for $p\to 0$ it is reduced to a Gustafson's
integral \cite{gus:some}). It was proved in \cite{rai:trans} after a
reduction to computations of fairly complicated sequences of determinants
on dense sets of parameters. Our proof is essentially shorter being based
on the completely elementary tools.
\end{remark}

\begin{remark}
As shown  in \cite{die-spi:elliptic}, integral
(\ref{ell-int-C}) implies validity of another highly non-trivial $C_n$
integral with six parameters (tagged as a type II integral), which generalizes
the Selberg integral and its $q$-analogues due to Gustafson \cite{gus:some}.
A direct derivation
of the latter integral was given in \cite{rai:trans}. It would be interesting
to consider it also from the present paper point of view.
\end{remark}

\section{A multiparameter $A_n$ integral}

For analyzing the type I $A_n$  elliptic beta integral
proposed in \cite{spi:elliptic}, we take $z=(z_1,\ldots,z_n)\in\C^n$, define $z_{n+1}$
via the relation $\prod_{i=1}^{n+1}z_i=1$, and introduce $2n+3$ complex parameters
$t=(t_1,\ldots,t_{n+1})$ and $s=(s_1,\ldots,s_{n+2})$.
The $A_n$ kernel function has the form
\begin{eqnarray}\nonumber
&& \rho(z,t,s;A_n)=
\prod_{i=1}^{n+1}\frac{\prod_{m=1}^{n+1}\Gamma(t_mz_i^{-1})
\prod_{j=1}^{n+2}\Gamma(s_jz_i)\, \Gamma(St_i)}
{\Gamma(TSz_i)\prod_{j=1}^{n+2}\Gamma(t_is_j)} % } &&
\\ && \makebox[4em]{}
\times \prod_{1\leq i<j\leq n+1}\frac{1}{\Gamma(z_iz_j^{-1},z_i^{-1}z_j)}
\frac{1}{\Gamma(T)}\prod_{j=1}^{n+2}\frac{\Gamma(STs_j^{-1})}{\Gamma(Ss_j^{-1})},
\label{kernel-A}\end{eqnarray}
where $T=\prod_{m=1}^{n+1}t_m$ and $S=\prod_{j=1}^{n+2}s_j$.
This function has poles at the points
$$
z_i=\{t_m q^jp^k,\, (TS)^{-1}q^{j+1}p^{k+1}\}, \; i=1,\ldots, n,\quad
z_{n+1}^{-1}=z_1\cdots z_n= \{s_lq^{j}p^{k}\},
$$
with $m=1,\ldots, n+1,\,l=1,\ldots, n+2,\, j,k \in\N,$ converging to zero and at
$$
z_i=\{s_l^{-1} q^{-j}p^{-k}\},\; i=1,\ldots, n,\quad
z_{n+1}^{-1}= \{t_m^{-1}q^{-j}p^{-k}, TSq^{-j-1}p^{-k-1}\},
$$
diverging to infinity.

\begin{theorem}
Suppose that $|t_m|, |s_l|<1$ and $|pq|<|TS|$. Then
\begin{equation}\label{ell-int-A}
\int_{\T^n}\rho(z,t,s;A_n)\frac{dz}{z}=\frac{(n+1)!(2\pi i)^n}
{(q;q)_\infty^n(p;p)_\infty^n}.
\end{equation}
\end{theorem}
\begin{proof}
The $A_n$-kernel function satisfies a $q$-difference equation, analogous
to (\ref{eqn}),
\begin{eqnarray}\nonumber
&& \rho(z,qt_1,t_2,\ldots,t_{n+1},s;A_n)-\rho(z,t,s;A_n)
\\ && \makebox[2em]{}
=\sum_{i=1}^n\left(g_i(z_1,\ldots,q^{-1}z_i,\ldots,z_n,t,s)-g_i(z,t,s)\right),
\label{eqn-A}\end{eqnarray}
where
\begin{equation}
\frac{g_i(z,t,s)}{\rho(z,t,s;A_n)}=
\prod_{j=1,\neq i}^{n+1}\frac{\theta(t_1z_j^{-1};p)}{\theta(z_iz_j^{-1};p)}
\prod_{j=1}^{n+2}\frac{\theta(z_is_j;p)}{\theta(t_1s_j;p)}
\frac{\theta(z_iT/t_1,TSt_1;p)}{\theta(T,TSz_i;p)}\frac{t_1}{z_i}.
\label{g-A}\end{equation}
Using the relation
\begin{eqnarray}\nonumber
&& \frac{\rho(\ldots,q^{-1}z_i,\ldots;A_n)}
{\rho(\ldots, z_i,\ldots;A_n)}=
\prod_{m=1}^{n+2}\frac{\theta(s_mz_{n+1};p)}{\theta(q^{-1}s_mz_i;p)}
\prod_{j=1}^{n+1}\frac{\theta(t_jz_i^{-1};p)}{\theta(q^{-1}t_jz_{n+1}^{-1};p)}
\frac{\theta(q^{-1}TSz_i;p)}{\theta(TSz_{n+1};p)}
\\ && \makebox[2em]{} \times
\prod_{j=1,\neq i}^n \frac{\theta(q^{-1}z_iz_j^{-1},q^{-1}z_jz_{n+1}^{-1};p)}
{\theta(z_i^{-1}z_j,z_j^{-1}z_{n+1};p)}
\frac{\theta(q^{-2}z_iz_{n+1}^{-1};p)z_i^2}
{\theta(z_iz_{n+1}^{-1};p)qz_{n+1}^2},
\label{eqn-A-zi}\end{eqnarray}
equation (\ref{eqn-A}) can be rewritten in the form
\begin{eqnarray}\nonumber
\lefteqn{\prod_{i=1}^{n+1}\frac{\theta(t_1z_i^{-1};p)}{\theta(TSz_i;p)}
\frac{\theta(t_1S;p)}{\theta(T;p)}\prod_{j=1}^{n+2}
\frac{\theta(TSs_j^{-1};p)}{\theta(t_1s_j;p)}-1 } &&
\\ \nonumber &&
=\frac{t_1\theta(t_1TS;p)}{\theta(T;p)\prod_{j=1}^{n+2}\theta(t_1s_j;p)}
\sum_{i=1}^n\frac{1}{z_i\theta(z_iz_{n+1}^{-1};p)}
\prod_{j=1,\neq i}^n\frac{\theta(t_1z_j^{-1};p)}{\theta(z_iz_j^{-1};p)}
\\ \nonumber && \makebox[2em]{} \times
\Biggl(\frac{z_i^{n+1}}{z_{n+1}^{n+1}}
\prod_{j=1}^{n+2}\theta(s_jz_{n+1};p)
\frac{\prod_{j=1}^{n+1}\theta(t_jz_i^{-1};p)\theta(q^{-1}z_iTt_1^{-1};p)}
{\prod_{j=2}^{n+1}\theta(q^{-1}t_jz_{n+1}^{-1};p)\theta(TSz_{n+1};p)}
\\ && \makebox[3em]{}
\times \prod_{j=1,\neq i}^n\frac{\theta(q^{-1}z_jz_{n+1}^{-1};p)}
{\theta(z_jz_{n+1}^{-1};p)} - \prod_{j=1}^{n+2}\theta(s_jz_i;p)
\frac{\theta(t_1z_{n+1}^{-1},z_iTt_1^{-1};p)}{\theta(TSz_i;p)} \Biggr).
\label{eqn-exp-A}\end{eqnarray}
Both sides of this equality are invariant under the transformation
$z_1\to pz_1$ and have equal sets of poles in $z_1$ and corresponding residues.
Indeed, it is not difficult to check coincidence of residues at
$z_1=(TS)^{-1}$ and we skip it. Comparison of the residues at
$z_{n+1}^{-1}=TS$ yields the equality
\begin{eqnarray}\nonumber
&& \sum_{i=1}^n
\theta(q^{-1}z_iTt_1^{-1};p)
\prod_{j=1,\neq i}^n\frac{\theta(TSz_jq^{-1};p)}
{\theta(z_i^{-1}z_j;p)}
\prod_{j=2}^{n+1}\theta(t_jz_i^{-1};p)
\\ &&\makebox[2em]{}
=\theta(t_1^{-1}S^{-1};p)\prod_{j=2}^{n+1}\theta(q^{-1}t_jTS;p).
\label{A-function}\end{eqnarray}
If we multiply both sides by
$\prod_{1\leq i<j\leq n } z_j\theta(z_i/z_j;p),$
then they become antisymmetric holomorphic theta functions
of $z_1,\dots,z_{n-1}$ obeying $A_{n-1}$ symmetry and  vanishing
at $z_{n}=z_i, i<n$. But any such function is proportional
to the given factor. Therefore the left-hand side of (\ref{A-function})
should be independent of $z_i$. Setting $z_1=q(TS)^{-1}$
we obtain the right-hand side expression.

Residues of the poles at $z_{n+1}=t_kq^{-1},\, k=2,\ldots,n+1,$
sum to the expression proportional to
\begin{equation}
\sum_{i=1}^n \frac{\prod_{j=2,\neq k}^{n+1}\theta(t_j/z_i;p)}
{\prod_{j=1,\neq i}^n \theta(z_j/z_i;p)}\theta(qt_1/z_iT;p),
\label{zn+1=tj/q}\end{equation}
which is equal to zero being the $TS=q/t_k$ subcase of (\ref{A-function}).

The poles at $z_1=z_j,1/z_j, j=2,\ldots, n,$ on the
right-hand side of (\ref{eqn-exp-A}) cancel each other.
Finally, comparison of the residues of the poles at $z_1=z_{n+1}$
yields the equality
\begin{equation}
\sum_{i=1}^n \theta(z_iT/qt_1;p)
\prod_{j=1,\neq i}^n\frac{\theta(z_j/qz_1;p)}{\theta(z_j/z_i;p)}
\prod_{j=2}^{n+1}\frac{\theta(t_j/z_i;p)}
{\theta(t_j/qz_1;p)}=\theta(z_1T/t_1;p).
\label{z1=zn+1}\end{equation}
After multiplication of both sides by
$\prod_{1\leq i<j\leq n } z_j \theta(z_i/z_j;p),$  they become
holomorphic $A_{n-2}$ antisymmetric functions of
$z_2,\ldots,z_{n-1}$ ($z_1$ is considered as a parameter)
vanishing at $z_j=z_1,z_n$ and therefore they are proportional to each
other up to a factor independent on $z_2,\dots,z_{n-1}$. For
$z_2=qz_1$, equality (\ref{z1=zn+1}) is true and, so, it is true in general.

From the listed properties it follows that the difference of two sides
of equality (\ref{eqn-exp-A}) does not depend on $z_1$. Setting $z_1=t_1$, we see that
this difference is equal to zero and (\ref{eqn-exp-A}) is true in general.

Integrating equation (\ref{eqn-A}) over the variables $z\in \T^n$, we obtain
\begin{equation}
I(qt_1,t_2,\ldots,t_{n+1},s)-I(t,s)=\sum_{i=1}^n
\left(\int_{\T^{i-1}\times(q^{-1}\T)\times \T^{n-i}}-\int_{\T^n}\right)
g_i(z,t,s)\frac{dz}{z},
\label{int-eqn-A}\end{equation}
where $I(t,s)=\int_{\T^n}\rho(z,t,s;A_n)dz/z$.

Poles of functions $g_i(z,t,s)$ are located at
\begin{eqnarray*}
&& z_i=\{ t_mq^jp^k,\; (TS)^{-1}q^{j}p^{k+1} \}, \quad i=1\ldots,n,\; m=1,\ldots,n+1,  \\
&& z_{n+1}^{-1}=\{ s_mq^jp^j \},\quad m=1,\ldots,n+2,
\end{eqnarray*}
converging to zero, and
\begin{eqnarray*}
&& z_i=\{ s_m^{-1}q^{-j-1}p^{-k} \},\quad i=1\ldots,n,\; m=1,\ldots,n+2, \\
&& z_{n+1}^{-1}=\{ t_1^{-1}q^{-j-1}p^{-k},\; t_m^{-1}q^{-j}p^{-k},\;
TSq^{-j-1}p^{-k-1} \},\quad m=2,\ldots,n+1,
\end{eqnarray*}
diverging to infinity for $j,k\in\mathbb{N}$.
Taking $|t_m|<|q|,\, m=2,\ldots,n+1$ and $|p|<|TS|$, we see that
the region $1\leq |z_i|\leq |q|^{-1}$ does not contain poles
and we can deform $q^{-1}\T$ to $\T$ in (\ref{int-eqn-A}) yielding
zero on the right-hand side. Thus, $I(qt_1,t_2,\ldots,t_{n+1},s)=I(t,s)$
in the taken parameter region.

Expansion of the integral in the infinite series in small $p$ allows the
iteration of the scaling $t_1\to pt_1$ termise until reaching the limiting point
$t_1=0$. Applying the analytical continuation procedure similar
to that used in the $C_n$ case, we see that $I(t,s)$ does not depend on $t_m$.
Integral (\ref{ell-int-A}) can be rewritten in the
$t_m\leftrightarrow s_m$ symmetric form
\begin{eqnarray}\nonumber
\lefteqn{ \int_{\T^n} \frac{\prod_{i=1}^{n+1}\prod_{m=1}^{n+2}
\Gamma(t_mz_i^{-1},s_mz_i)}
{\prod_{1\leq i<j\leq n+1}\Gamma(z_iz_j^{-1},z_i^{-1}z_j)}\frac{dz}{z}
} && \\ &&
= \frac{(n+1)!(2\pi i)^n}{(q;q)_\infty^n(p;p)_\infty^n}
\prod_{m=1}^{n+2}\Gamma(Ss_m^{-1},At_m^{-1})
\prod_{j,k=1}^{n+2}\Gamma(t_js_k).
\label{ell-int-A-sym}\end{eqnarray}
where $S=\prod_{m=1}^{n+2}s_m, A=\prod_{m=1}^{n+2}t_m,
AS=pq$. Therefore $I(t,s)$ does not depend on $s_m$ as well. The explicit
form of $I=I(q,p)$, given by the right-hand side of (\ref{ell-int-A}),
is found from the residue calculus described in \cite{spi:elliptic}.
\end{proof}

\begin{remark}
Formula (\ref{ell-int-A}) was conjectured in \cite{spi:elliptic}
with various justifications. For $p\to 0$ it is reduced to a Gustafson's
integral \cite{gus:some} and its first proof was obtained in
\cite{rai:trans} with the determinantal technique similar to that
used for treating the multiparameter $C_n$ integral.
Again, our proof looks simpler.
\end{remark}

\begin{remark}
As shown in \cite{spi:elliptic}, a combination of (\ref{ell-int-C}) and (\ref{ell-int-A})
integrals implies validity of two type II exact integration formulas
(which are conventionally considered as the $A_n$ root system objects), having
different expressions for odd and even multiplicities of integrations. Their
direct analysis from the present point of view would be of interest as well.
\end{remark}

\section{Modified integrals with a base on the unit circle}

All the elliptic beta integrals discussed so far are defined only
when both bases lie within the unit circle $|q|,|p|<1$. Elliptic
hypergeometric integrals which are well defined when one of the
bases lies on the unit circle, say $|q|=1$, were introduced in
\cite{spi:elliptic}. They are constructed with the help
of the modified elliptic gamma function (\ref{ell-d}), (\ref{gamma-tr}).
A unit circle analog of the type II six parameter $C_n$ elliptic beta
integral was constructed in \cite{die-spi:elliptic}. Below we describe
several new unit circle elliptic beta integrals of type I for the root
systems $A_n$ and $C_n$.

\begin{theorem}
For the $C_n$ root system we take $n$ variables $u=(u_1,\ldots,u_n)$ $\in\C^n$,
$2n+3$ complex parameters $g=(g_1,\ldots,g_{2n+3})$, and define the kernel
function
\begin{eqnarray}\nonumber
&& \rho(u,g;C_n)=
\prod_{i=1}^n\frac{\prod_{m=1}^{2n+3}G(g_m\pm u_i)}
{G(\pm 2u_i,\mathcal{A}\pm u_i)}\frac{\prod_{m=1}^{2n+3}
G(\mathcal{A}-g_m)}{\prod_{1\leq m<s\leq 2n+3}G(g_m+g_s)}
\\ && \makebox[6em]{} \times
\prod_{1\leq i<j\leq n}\frac{1}{G(\pm u_i\pm u_j)},
\label{kernel-C-unit}\end{eqnarray}
where $\mathcal{A}=\sum_{m=1}^{2n+3}g_m$, the modified elliptic gamma function
 $G(u)\equiv G(u;\boldsymbol{\omega})$, and
$$
G(c\pm a \pm b)\equiv G(c+a+b,c+a-b,c-a+b,c-a-b).
$$
Now we suppose that $\text{Im} (\omega_1/\omega_2)\geq 0$ and $\text{Im}
(\omega_3/\omega_1), \text{Im} (\omega_3/\omega_2) >0$, and
$$
\text{Im}(g_m/\omega_3)<0,\; m=1,\ldots,2n+3, \qquad
\text{Im}((\mathcal{A}-\omega_1-\omega_2)/\omega_3)>0.
$$
Then
\begin{equation}
\int_{-\frac{\omega_3}{2}}^{\frac{\omega_3}{2}}\cdots
\int_{-\frac{\omega_3}{2}}^{\frac{\omega_3}{2}} \rho(u,g;C_n)
\frac{du_1}{\omega_2}\cdots \frac{du_N}{\omega_2}
= \kappa^n 2^n n!
\label{C_n-circle}\end{equation}
with
\begin{equation}
\kappa= \frac{-(\tilde q;\tilde q)_\infty}
{(q;q)_\infty(p;p)_\infty(r;r)_\infty}
\label{kappa}\end{equation}
and the integration along the straight line connecting points $-\omega_3/2$
and $\omega_3/2$.
\end{theorem}
\begin{proof}
After substituting relation (\ref{gamma-tr}) in (\ref{kernel-C-unit}) and some
straightforward computations, we obtain
\begin{eqnarray}\nonumber
&& \rho(u,g;C_n)=e^{-\pi i nP(0)}
\prod_{1\leq i<j\leq n}\frac{1}{\Gamma(z_i^\pm z_j^\pm;\tilde r,\tilde p)}
\prod_{i=1}^n\frac{\prod_{m=1}^{2n+3}\Gamma(t_mz_i^\pm;\tilde r,\tilde p)}
{\Gamma(z_i^{\pm 2},Az_i^\pm;\tilde r,\tilde p)}
\\ && \makebox[4em]{} \times
\frac{\prod_{m=1}^{2n+3}\Gamma(At_m^{-1};\tilde r,\tilde p)}
{\prod_{1\leq m<s\leq 2n+3}\Gamma(t_mt_s;\tilde r,\tilde p)},
\label{rho-C'}\end{eqnarray}
where $z_j=e^{-2\pi i u_j/\omega_3},$ $t_m=e^{-2\pi i g_m/\omega_3}$, and
$A=\prod_{m=1}^{2n+3}t_m$. Since
$$
du_j=\frac{\omega_3}{2\pi i}\frac{dz_j}{z_j},
$$
the integrals in $u_j$ become equivalent to integrals in $z_j$ over the
unit circle $\T$. Applying $C_n$ integration formula (\ref{ell-int-C})
with $C=\T$ (it is allowed due to the constraints $|t_m|<1, |pq|<|A|$
which are assumed to hold) and $q,p$ replaced by $\tilde r,\tilde p$,
we obtain that the left-hand side of (\ref{C_n-circle}) is equal to
$\kappa^n 2^n n!,$ where
$$
\kappa=\frac{\omega_3e^{\frac{\pi i}{12}(\sum_{m=1}^3\omega_m)
(\sum_{m=1}^3\omega_m^{-1}) }}
{\omega_2(\tilde r; \tilde r)_\infty(\tilde p;\tilde p)_\infty}.
$$
Using the transformation property of the Dedekind $\eta$-function
$$
e^{-\frac{\pi i}{12\tau}} \left(e^{-2\pi i/\tau};e^{-2\pi
i/\tau}\right)_\infty =(-i\tau)^{1/2}e^{\frac{\pi i\tau}{12}}
\left(e^{2\pi i\tau};e^{2\pi i\tau}\right)_\infty ,
$$
it can be shown \cite{die-spi:elliptic} that $\kappa$ takes the
form presented in (\ref{kappa}).
\end{proof}

\begin{theorem}
For the simplest $A_n$ root system integral, we take integration variables
$u=(u_1,\ldots,u_n)$ $\in\C^n$
and complex parameters $g=(g_1,\ldots,g_{n+1})$, $h=(h_1,\ldots,h_{n+2})$.
Denoting $u_{n+1}=-\sum_{k=1}^n u_k$, we define the kernel function
\begin{eqnarray}\nonumber
&& \rho(u,g,h;A_n)=
\prod_{i=1}^{n+1}\frac{\prod_{m=1}^{n+1}G(g_m-u_i)
\prod_{j=1}^{n+2}G(h_j+u_i)\, G(H+g_i)}
{G(F+H+u_i)\prod_{j=1}^{n+2}G(g_i+h_j)}
\\ && \makebox[4em]{}
\times \prod_{1\leq i<j\leq n+1}\frac{1}{G(u_i-u_j,u_j-u_i)}
\frac{1}{G(F)}\prod_{j=1}^{n+2}\frac{G(F+H-h_j)}{G(H-h_j)},
\label{kernel-A-unit}\end{eqnarray}
where $F=\sum_{m=1}^{n+1}g_m$ and $H=\sum_{j=1}^{n+2}h_j$.
Supposing that
\begin{eqnarray*}
&& \text{Im}\, \frac{\omega_1}{\omega_2}\geq 0,\quad
\text{Im}\, \frac{\omega_3}{\omega_1} >0,\quad
\text{Im}\, \frac{\omega_3}{\omega_2} >0, \\
&& \text{Im}\, \frac{g_m}{\omega_3}<0,\quad
\text{Im}\, \frac{h_j}{\omega_3}<0,\quad
\text{Im}\, \frac{F+H-\omega_1-\omega_2}{\omega_3}>0
\end{eqnarray*}
for $m=1,\ldots,n+1,\, j=1,\ldots,n+2$, we have
\begin{equation}
\int_{-\frac{\omega_3}{2}}^{\frac{\omega_3}{2}}\cdots
\int_{-\frac{\omega_3}{2}}^{\frac{\omega_3}{2}} \rho(u,g,h;A_n)
\frac{du_1}{\omega_2}\cdots \frac{du_N}{\omega_2}
= \kappa^n (n+1)!
\label{A_n-circle}\end{equation}
with the same $\kappa$ as in the previous theorem.
\end{theorem}
\begin{proof}
Substituting relation (\ref{gamma-tr}) in (\ref{kernel-A-unit}), we obtain
\begin{eqnarray}\nonumber
e^{\pi i nP(0)}\lefteqn{ \rho(u,g,h;A_n)=
\prod_{i=1}^{n+1}\frac{\prod_{m=1}^{n+1}\Gamma(t_mz_i^{-1};\tilde r,\tilde p)
\prod_{j=1}^{n+2}\Gamma(s_jz_i;\tilde r,\tilde p)\, \Gamma(St_i;\tilde r,\tilde p)}
{\Gamma(TSz_i;\tilde r,\tilde p)\prod_{j=1}^{n+2}\Gamma(t_is_j;\tilde r,\tilde p)}
}&& \\ &&
\times \prod_{1\leq i<j\leq n+1}\frac{1}{\Gamma(z_iz_j^{-1},z_jz_i^{-1};\tilde r,\tilde p)}
\frac{1}{\Gamma(T;\tilde r,\tilde p)}\prod_{j=1}^{n+2}
\frac{\Gamma(TSs_j^{-1};\tilde r,\tilde p)}{\Gamma(Ss_j^{-1};\tilde r,\tilde p)},
\label{rho-A'}\end{eqnarray}
where $z_j=e^{-2\pi i u_j/\omega_3},$ $t_m=e^{-2\pi i g_m/\omega_3},$
$s_j=e^{-2\pi i h_j/\omega_3},$ $T=\prod_{m=1}^{n+1}t_m,$
$S=\prod_{j=1}^{n+2}s_j.$
Similar to the $C_n$ case, the integrals in $u_j$ become equivalent to integrals
in $z_j$ over the unit circle $\T$. Applying $A_n$ integration formula (\ref{ell-int-A})
(which is allowed since our constraints are equivalent to
$|t_m|<1,$ $|s_l|<1,$ and $|pq|<|TS|$) with $q,p$ replaced by $\tilde r,\tilde p$,
we obtain that the integral in (\ref{A_n-circle}) equals to $\kappa^n (n+1)!.$
\end{proof}

\begin{theorem}
For the unit circle analog of the type I elliptic beta integral introduced in
\cite{spi-war:inversions}, we take integration variables $u=(u_1,\ldots,u_n)$ $\in\C^n$
and complex parameters $g=(g_1,\ldots,g_{n+3})$, $h=(h_1,\ldots,h_n)$.
Denoting $u_{n+1}=-\sum_{k=1}^n u_k$, we define the kernel function
\begin{eqnarray}\nonumber
\lefteqn{ \delta(u,g,h;A_n)=\prod_{1\leq i<j\leq n+1}\frac{G(F-u_i-u_j)}
{G(u_i-u_j,u_j-u_i)}
\prod_{j=1}^{n+3}\prod_{m=1}^n\frac{G(F+h_m-g_j)}{G(h_m+g_j)}
} && \\ && \times
\prod_{i=1}^{n+1}\frac{\prod_{m=1}^n G(h_m+u_i)
\prod_{j=1}^{n+3}G(g_j-u_i)}{\prod_{m=1}^n G(F+h_m-u_i)}
\prod_{1\leq j<k\leq n+3}
\frac{1}{G(F-g_j-g_k)},
\label{kernel-D-unit}\end{eqnarray}
where $F=\sum_{j=1}^{n+3}g_j$.
Supposing that $\text{Im} (\omega_1/\omega_2)\geq 0$ and $\text{Im}
(\omega_3/\omega_1),$ $\text{Im} (\omega_3/\omega_2) >0$, and
$\text{Im}(g_j/\omega_3)<0,\;$ $\text{Im}(h_m/\omega_3)<0,\;$
$\text{Im}((F+h_m-\omega_1-\omega_2)/\omega_3)>0$
for $j=1,\ldots,n+3,\, m=1,\ldots,n,$ we have
\begin{equation}
\int_{-\frac{\omega_3}{2}}^{\frac{\omega_3}{2}}\cdots
\int_{-\frac{\omega_3}{2}}^{\frac{\omega_3}{2}} \delta(u,g,h;A_n)
\frac{du_1}{\omega_2}\cdots \frac{du_N}{\omega_2}
= \kappa^n (n+1)!
\label{D_n-circle}\end{equation}
with the same $\kappa$ as in the previous cases.
\end{theorem}
\begin{proof}
Substituting relation (\ref{gamma-tr}) in (\ref{kernel-D-unit}), we obtain
$$
\delta(u,g,h;A_n)=e^{-\pi i nP(0)}\Delta(z,t,s;A_n),
$$
where
\begin{eqnarray}\label{rho-D'}
\lefteqn{ \Delta(z,t,s;A_n)=\prod_{1\leq i<j\leq n+1}
\frac{\Gamma(Dz_i^{-1}z_j^{-1};\tilde r,\tilde p)}
{\Gamma(z_iz_j^{-1},z_i^{-1}z_j;\tilde r,\tilde p)}
\prod_{j=1}^{n+3}\prod_{m=1}^n\frac{\Gamma(Dt_ms_j^{-1};\tilde r,\tilde p)}
{\Gamma(t_ms_j;\tilde r,\tilde p)}
 }&&  \\ &&
\times\prod_{i=1}^{n+1}\frac{\prod_{m=1}^{n}\Gamma(t_mz_i;\tilde r,\tilde p)
\prod_{j=1}^{n+3}\Gamma(s_jz_i^{-1};\tilde r,\tilde p)}
{\prod_{m=1}^{n}\Gamma(Dt_mz_i^{-1};\tilde r,\tilde p)}
\prod_{1\leq j<k\leq n+3}
\frac{1}{\Gamma(Ds_j^{-1}s_k^{-1};\tilde r,\tilde p)}
\nonumber\end{eqnarray}
with $z_j=e^{-2\pi i u_j/\omega_3},$ $s_j=e^{-2\pi i g_j/\omega_3},$
$t_m=e^{-2\pi i h_m/\omega_3},$ $D=\prod_{j=1}^{n+3}s_j.$
As shown by Warnaar and the author in \cite{spi-war:inversions},
for $|s_j|, |t_m|<1$ and $|\tilde p\tilde r|<|Dt_m|$ (which is true by the taken
assumptions) the following $A_n$ (or $``D_n"$) elliptic beta integral holds
\begin{equation}
\int_{\T^n}\Delta(z,t,s;A_n)\frac{dz}{z}
=\frac{(2\pi i)^n (n+1)!}{(\tilde r;\tilde r)_\infty^n
(\tilde p;\tilde p)_\infty^n}.
\label{ell-int-D}\end{equation}
Applying this formula in our case, we obtain that the left-hand side
of (\ref{D_n-circle}) equals to $\kappa^n(n+1)!$ as required.
\end{proof}

\section{$q$-Reductions of the modified integrals}

In this section we consider $q$-reductions of the modified elliptic beta integrals
described in the previous section. They correspond to the limit
$\text{Im}(\omega_3)\to \infty$ taken in such a way that $p, r\to 0$.
In this limit the modified elliptic gamma functions
reduce to the double sine functions and we obtain some beta type integrals
over infinite contours (for earlier examples of non-compact $q$-beta
integrals of such type, see, e.g., \cite{die-spi:elliptic} and references therein).
However, this is a formal limit and the resulting
integrals need to be proved rigorously.

\begin{theorem}
For the $C_n$ root system we take the notation used in formula
(\ref{kernel-C-unit}) and define the kernel function
\begin{eqnarray}\nonumber
&& \rho(u,g;\omega;C_n)=\prod_{1\leq i<j\leq n}S(\pm u_i\pm u_j)
\prod_{i=1}^n\frac{S(\pm 2u_i,\mathcal{A}\pm u_i)}
{\prod_{m=1}^{2n+3}S(g_m\pm u_i)}
\\ &&\makebox[4em]{}\times
\frac{\prod_{1\leq m<s\leq 2n+3}S(g_m+g_s)}
{\prod_{m=1}^{2n+3}S(\mathcal{A}-g_m)},
\label{kernel-C-q}\end{eqnarray}
where $\mathcal{A}=\sum_{m=1}^{2n+3}g_m$ and $S(u)\equiv S(u;\boldsymbol{\omega})$
is the double sine function. Now we suppose that
$\text{Im} (\omega_1/\omega_2)\geq 0$ and $\text{Re} (\omega_1/\omega_2)>0$
together with
$$
\text{Re} (g_m/\omega_2) > 0, \quad
\text{Re}((\mathcal{A}-\omega_1)/\omega_2)<1
$$
for $m=1,\dots,2n+3.$ Then
\begin{equation}
\int_{\L^n}\rho(u,g;\omega;C_n)
\frac{du_1}{\omega_2}\cdots \frac{du_N}{\omega_2}
= (-2)^n n!\frac{(\tilde q;\tilde q)_\infty^n}{(q;q)_\infty^n}
\label{C_n-q}\end{equation}
with the integration contour $\L=i\omega_2\R$.
\end{theorem}
\begin{proof}
First, we consider convergence of the taken integral.
For this we assume that the integration domain is bounded to an $n$-dimensional
cube and consider integration in the variable $u_k=i\omega_2x_k$ for some
fixed $k$ when the size of the cube goes to infinity. Due to
the $u_j\to-u_j$ symmetry, it is sufficient to consider
the subdomain $0\leq x_j<\infty,\, j=1,\dots,n$.  We impose
the constraint $\text{Re} (\omega_1/\omega_2)>0$,
so that for $x_k\to+\infty$ we have $\text{Re} (u_k/i\omega_1)\to+\infty$.
From the asymptotic behavior of $S(u)$ described in (\ref{asymp1}), (\ref{asymp2}),
we obtain
$$
\frac{S(\pm 2u_k,\mathcal{A}\pm u_k)}
{\prod_{m=1}^{2n+3}S(g_m\pm u_k)}=O(e^{-2(n-1)\pi ix_k^2\omega_2/\omega_1
-2\pi n x_k(1+\omega_2/\omega_1)}).
$$
Using the reflection formula
$$
S(u)S(-u)=e^{-\pi iB_{2,2}(u)}(1-e^{-2\pi iu/\omega_2})
(1-e^{-2\pi iu/\omega_1}),
$$
for $x_k\to+\infty$ we find the estimate
$$
\prod_{1\leq i<j\leq n}S(\pm u_i\pm u_j)=
O(e^{-2\pi i(n-1)\sum_{j=1}^nB_{2,2}(-u_j)})
$$
independently on the values of other $x_j$. Thus, on the
boundary of the integration domain we have
$$
\rho(u,g;\omega;C_n)=O(e^{-2\pi (1+\omega_2/\omega_1)\sum_{j=1}^nx_j}),
$$
i.e., the kernel function decays exponentially fast and the integral converges.

After taking the limit $p\to 0$ in (\ref{eqn-C}), we obtain the
difference equation for our kernel function (scalings of $t_m$ and $z_i$
by $q^{\pm1}$ are equivalent to the shifts of $g_m$ and $u_i$
by $\pm\omega_1$):
\begin{eqnarray} \nonumber
&& \rho(u,g_1+\omega_1,g_2,\ldots,g_{2n+3};\omega;C_n)-
\rho(u,g;\omega;C_n)
\\ && \makebox[4em]{}
=\sum_{i=1}^n\left(f_i(u_1,...,u_i-\omega_1,\ldots,u_n,g)-f_i(u,g)\right),
\label{eqn-C-q}\end{eqnarray}
where
\begin{eqnarray}\nonumber
&& f_i(u,g)=\rho(u,g;\omega;C_n)\prod_{j=1,\neq i}^n\frac{(1-t_1z_j)(1-t_1/z_j)}
{(1-z_iz_j)(1-z_i/z_j)}
\\ && \makebox[4em]{} \times
\frac{\prod_{j=1}^{2n+3}(1-t_jz_i)}{\prod_{j=2}^{2n+3}(1-t_1t_j)}
\frac{(1-t_1A)t_1}{(1-z_i^2)(1-Az_i)z_i}
\label{g-C-q}\end{eqnarray}
with $z_j=e^{2\pi iu_j/\omega_2},\, t_k=e^{2\pi ig_k/\omega_2}$,
$A=e^{2\pi i\mathcal{A}/\omega_2}.$ Since $S(u)$ is symmetric in $\omega_{1,2}$,
we have a similar equation for shifts by $\omega_2$.

We suppose now temporarily that quasiperiods $\omega_{1,2}$ are
real, incommensurate, and $\omega_1<\omega_2$. We impose also the
constraint $\text{Re}(\mathcal{A})<\omega_1$.
Integrating equation (\ref{eqn-C-q}) over $u_j$ along the contour $\L$, we obtain
$$
I(g_1+\omega_1,g_2,\dots,g_{2n+3})
=I(g)\equiv\int_\L \rho(u,g;\omega;C_n)du_1\cdots du_n.
$$
This follows from
the observation that under taken constraints the functions $f_i(u,g)$ do not
have poles in the strip $-\omega_1<u_i<0$ and we can shift the contour
$\L-\omega_1$ to $\L$. By symmetry, we have
$I(g_1+\omega_2,g_2,\dots,g_{2n+3})=I(g).$ After an appropriate finite
deformation of the integration contour we can reach the equality
$I(g_1+j\omega_1-k\omega_2,g_2,\dots,g_{2n+3})=I(g)$ for all $j,k\in\N$
such that $j\omega_1-k\omega_2\in[0,\omega_2]$. Due to the denseness
of this set of points, we conclude that the function $I(g)$ does not depend
on $g_1$ and, so, on all $g_m$. As a result, its explicit form,
given by the left-hand side of (\ref{C_n-q}), can be found from the
residue calculus analogous to those described in \cite{die-spi:elliptic}.
Finally, by analyticity we can extend permitted values of parameters
to the domain indicated in the formulation of the theorem.
\end{proof}

\begin{theorem}
For the simplest $A_n$ root system integral, we take the notation for parameters
as in formula (\ref{kernel-A-unit}) and define the kernel function
\begin{eqnarray}\nonumber
&& \rho(u,g,h;\omega;A_n)=
\prod_{i=1}^{n+1}\frac{S(F+H+u_i)\prod_{j=1}^{n+2}S(g_i+h_j)}
{\prod_{m=1}^{n+1}S(g_m-u_i)\prod_{j=1}^{n+2}S(h_j+u_i)\, S(H+g_i)}
\\ && \makebox[4em]{}
\times \prod_{1\leq i<j\leq n+1}S(u_i-u_j,u_j-u_i)\;
S(F)\prod_{j=1}^{n+2}\frac{S(H-h_j)}{S(F+H-h_j)}.
\label{kernel-A-q}\end{eqnarray}
Now we suppose that $\text{Im} (\omega_1/\omega_2)\geq 0$ and
$\text{Re} (\omega_1/\omega_2)>0$ together with
$$
\text{Re} (g_m/\omega_2) > 0, \quad \text{Re} (h_j/\omega_2) > 0,\quad
\text{Re}((F+H-\omega_1)/\omega_2)<1.
$$
Then
\begin{equation}
\int_{\L^n}\rho(u,g,h;\omega;A_n)
\frac{du_1}{\omega_2}\cdots \frac{du_N}{\omega_2}
= (-1)^n (n+1)!\frac{(\tilde q;\tilde q)_\infty^n}{(q;q)_\infty^n}
\label{A_n-q}\end{equation}
with the integration contour $\L=i\omega_2\R$.
\end{theorem}
\begin{proof}
Convergence of the integral under consideration follows in the way
similar to the previous case.

Needed difference equation for the kernel function is obtained from  the
$p\to 0$ limit of (\ref{eqn-A}):
\begin{eqnarray}\nonumber
&& \rho(u,g_1+\omega_1,g_2,\ldots,g_{n+1},h;\omega;A_n)-\rho(u,g,h;\omega;A_n)
\\ && \makebox[2em]{}
=\sum_{i=1}^n\left(f_i(u_1,\ldots,u_i-\omega_1,\ldots,u_n,g,h)-f_i(u,g,h)\right),
\label{eqn-A-q}\end{eqnarray}
where
\begin{equation}
\frac{f_i(u,g,h)}{\rho(u,g,h;\omega;A_n)}=
\prod_{j=1,\neq i}^{n+1}\frac{1-t_1/z_j}{1-z_i/z_j}
\prod_{j=1}^{n+2}\frac{1-z_is_j}{1-t_1s_j}
\frac{(1-z_iT/t_1)(1-TSt_1)t_1}{(1-T)(1-TSz_i)z_i}
\label{g-A-q}\end{equation}
with
$z_j=e^{2\pi iu_j/\omega_2},\, t_k=e^{2\pi ig_k/\omega_2}$,
$s_j=e^{2\pi ih_j/\omega_2}$, and
$A=e^{2\pi i\mathcal{A}/\omega_2}.$ By symmetry,
we have also a similar equation for the shifts by $\omega_2$.

We suppose now temporarily that quasiperiods $\omega_{1,2}$ are
positive (say, $\omega_1<\omega_2$) and incommensurate. We impose also the
constraints $\text{Re}(F+H)<\omega_2$ and $\text{Re}(g_m)>\omega_1,\,
m=2,\dots,n+1$.
Integrating equation (\ref{eqn-A-q}) over $u_j$ along the contour $\L$, we
obtain $I(g_1+\omega_1,g_2,\dots,g_{2n+3},h)$
$=I(g,h)\equiv\int_\L \rho(u,g,h;\omega;A_n)du_1\cdots du_n$ (the functions
$f_i(u,g,h)$ do not have poles in the strip $-\omega_1<u_i<0$).

Similar to the $C_n$ case, for an appropriate choice of the
integration contours we have
$I(g_1+j\omega_1-k\omega_2,g_2,\dots,g_{2n+3},h)=I(g,h)$ for all $j,k\in\N$
such that $j\omega_1-k\omega_2\in[0,\omega_2]$. As a result, the function
$I(g,h)$ does not depend on $g_1$ and, so, on all $g_m$.
Applying the residue calculus, analogous to
that described in \cite{spi:elliptic}, we find that $I(g,h)$ does not
depend on $h$ as well and equals to the expression on
the right-hand side of (\ref{A_n-q}). Again, analytical continuation
in parameters from restricted values to the admissible domain
completes the proof of the integral in question.
\end{proof}

\begin{remark}
The Gustafson's plain hypergeometric $A_n$ integral \cite{gus:some} is relevant for
the ordinary quantum Toda chain model (an observation due to S. Kharchev). Therefore,
it is natural to expect that a simplified version of the integral (\ref{A_n-q})
is related to normalization conditions of the Hamiltonian eigenfunctions
for the quantum multiparticle $q$-Toda chain of \cite{kls:unitary}.
\end{remark}

We shall not consider here unit circle analogs and relevant $q$-degenerations for
two other types of the $A_n$ integrals introduced in \cite{spi:elliptic}.

\begin{theorem}
For derivation of the $q$-reduction of integral (\ref{D_n-circle}), we take the
notation for parameters as in (\ref{kernel-D-unit}) and define the kernel function
\begin{eqnarray}\nonumber
\lefteqn{ \delta(u,g,h;\omega;A_n)=
\prod_{1\leq i<j\leq n+1}\frac{S(u_i-u_j,u_j-u_i)}{S(F-u_i-u_j)}
\prod_{j=1}^{n+3}\prod_{m=1}^n\frac{S(h_m+g_j)}{S(F+h_m-g_j)}
}&& \\ && \times
\prod_{i=1}^{n+1}\frac{\prod_{m=1}^n S(F+h_m-u_i)}
{\prod_{m=1}^n S(h_m+u_i)\prod_{j=1}^{n+3}S(g_j-u_i)}
\prod_{1\leq j<k\leq n+3}S(F-g_j-g_k).
\label{kernel-D-q}\end{eqnarray}
Now we suppose that $\text{Im} (\omega_1/\omega_2)\geq 0$ and
$\text{Re} (\omega_1/\omega_2)>0$ together with
$$
\text{Re} (g_j/\omega_2) > 0, \quad \text{Re} (h_m/\omega_2) > 0, \quad
\text{Re}((F+h_m-\omega_1)/\omega_2)<1.
$$
Then
\begin{equation}
\int_{\L^n}\delta(u,g,h;\omega;A_n)
\frac{du_1}{\omega_2}\cdots \frac{du_N}{\omega_2}
= (-1)^n (n+1)!\frac{(\tilde q;\tilde q)_\infty^n}{(q;q)_\infty^n}
\label{D_n-q}\end{equation}
with the integration contour $\L=i\omega_2\R$.
\end{theorem}

Convergence of the integral (\ref{D_n-q}) is deduced analogously to the previous
cases. Further steps of the proof of this
theorem are similar to the described $A_n$ and $C_n$ situations. Using
the equation for kernel function appearing after taking the $p\to 0$ limit
in the equation constructed in \cite{spi-war:inversions}, it can be
shown that the left-hand side of (\ref{D_n-q}) is not changed after the
shifts $g_m\to g_m+\omega_{1,2}$. Such shifts can be iterated after an appropriate
deformation of the contour of integration and this leads to the independence of
the integral on parameters $g_m$. Application of the appropriate residue calculus
(the $p\to 0$ partner of that used in \cite{spi-war:inversions}) leads
to equality (\ref{D_n-q}). We skip the details for they are simple enough and repetitive
to those given above for other root systems.

\medskip
\centerline{Acknowledgments}
\smallskip

This paper was written under the influence of Richard Askey's thread to
search for elementary treatment of special functions.
The author is indebted to J. F. van Diejen, S. M. Kharchev, E. M. Rains, and
S. O. Warnaar for stimulating discussions.

\end{document}